\documentclass{article}
\usepackage{arxiv}
\usepackage{amsfonts,amssymb,amsmath,amsthm}
\usepackage{graphicx,bm}
\usepackage{algorithm,algorithmic}
\usepackage{epstopdf}
\usepackage[labelformat=simple]{subfig}

\graphicspath{ {./figs/} }
\usepackage[labelformat=simple]{subfig}

\graphicspath{ {./figs/} }

\newtheorem{thm}{Theorem}
\newtheorem{remark}[thm]{Remark}
\newcommand {\av}[1] {\mbox{$\left\{\!\!\left\{ #1 \right\}\!\!\right\}$}}
\newcommand {\jump}[1] {\mbox{$\left[\!\left[ #1 \right]\!\right]$}}
\usepackage{doi}

\usepackage{hyperref}
\makeindex
  \title{Reduced-Order Modeling for Heston Stochastic Volatility Model}

\author{
\href{https://orcid.org/0000-0002-8136-0328}{\includegraphics[scale=0.06]{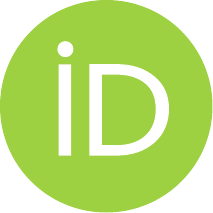}\hspace{1mm}Sinem Kozp{\i}nar}\\
   Department of Insurance, Ba\c{s}kent University\\
     Ankara-Turkey\\ \
     \texttt{sinemkozpinar@ankara.baskent.edu.tr}\\
          \And
  \href{https://orcid.org/0000-0001-5262-063X}{\includegraphics[scale=0.06]{orcid.pdf}\hspace{1mm}Murat Uzunca} \\
   Department of Mathematics, Sinop University\\
     Sinop-Turkey\\ \
     \texttt{muzunca@sinop.edu.tr}\\
     \And
     \href{https://orcid.org/0000-0003-1037-5431}{\includegraphics[scale=0.06]{orcid.pdf}\hspace{1mm}B\"ulent Karas\"ozen} \\
     Institute of Applied Mathematics \& Department of Mathematics\\
     Middle East Technical University\\
     Ankara-Turkey\\
     \texttt{bulent@metu.edu.tr}
}
\begin{document}

\maketitle

\begin{abstract}
In this paper, we compare the intrusive proper orthogonal decomposition (POD) with Galerkin projection and the data-driven dynamic mode decomposition (DMD), for Heston's option pricing model. The full order model is obtained by discontinuous Galerkin discretization in space and backward Euler in time. Numerical results for butterfly spread, European and digital call options reveal that in general DMD requires more modes than the POD modes for the same level of accuracy. However, the speed-up factors are much higher for DMD than POD due to the non-intrusive nature of the DMD.	
\end{abstract}

\keywords{
Option pricing, Heston model, discontinuous Galerkin method, proper orthogonal decomposition, reduced order modeling, dynamic mode decomposition\\
MSC Classification 2000. 65M60,  91B25,  91G80}

\section{Introduction}

Under Heston framework \cite{Heston01041993}, the volatility is treated as a square root process which enables us to predict the volatility smirk arising in the option prices. As a result, one can obtain more accurate prices for European currency options.  One of the well-known techniques to solve option pricing problems is the Monte Carlo integration.
Although its implementation is straightforward, it requires a large number of realizations to achieve high accuracy. Moreover, one can estimate only one option price for a given underlying value. Therefore, the simulations can be too costly which leads us to the discretization methods based on the solution of partial differential equations (PDEs). Option pricing problems under the Heston model are also represented as diffusion-convection-reaction PDEs \cite{Heston01041993, Lipton01mathematical}.
The diffusion matrix and convective field depend on the volatility.  The diffusion matrix contains cross-diffusion terms as a result of the correlation between the volatility and the underlying security price. In addition, the initial and boundary data are discontinuous and less regular for different options.
Heston's PDE is solved in the literature by finite differences \cite{Hout10adi,Tangman08,During14} and by finite elements \cite{Winkler01}. Recently  Heston's PDE has been  investigated for various option pricing models applying  discontinuous Galerkin (dG) method by the authors \cite{Kozpinar20}. The discontinuous Galerkin finite element (dGFEM) method has several advantages over finite difference  and continuous finite element  methods for solving diffusion-convection-reaction equations.
Because the dGFEM does not require continuity across the inter-element boundaries, the number of degrees of freedom is larger than the continuous FEMs. On the other hand, the dGFEM has a number of desirable properties like the treatment of the convective term by upwinding, weakly enforcement of the boundary conditions, and ease of parallelization.

In the last two decades, reduced-order modeling became an important tool for simulating engineering problems efficiently. 
Reduced-order models (ROMs) reduce the computational complexity and time by approximating the full order model (FOM) of the high dimensional discretized PDEs as lower dimensional models. 
This enables fast simulation based studies like calibration and hedging \cite{Cont11,Pironneau12pricing,Sachs13} in option pricing. 
In this paper, we compare two model order reduction (MOR) methods, the proper orthogonal decomposition (POD) and dynamic mode decomposition (DMD) for Heston's option pricing model.
The POD with Galerkin projection is the most known and frequently used reduced-order modeling method. 
POD based reduced-order modeling is investigated for pricing European option and American option wherein different from European option there exist variational inequalities \cite{cen22ipg,zeng21nqh,zeng21wpo} to be handled, under Black-Scholes and Heston model \cite{Balajewicz16,Burkovska15,Mayerhofer14,Peherstorfer15reduced}.  
Reduced-order models based on  POD are optimal in terms of energy content. The energy content is important but is not sufficient in general to catch the dynamical behavior. 
In this paper, we apply the data-driven reduced-order modeling technique DMD \cite{Schmid10,Tu14} for option pricing and compare it with the POD with respect to accuracy and speed-up over the full order models. The DMD is able to extract dynamically relevant flow features from time-resolved experimental or numerical data by generalizing the global stability modes and approximating the eigenvalues and eigenfunctions of the Koopman operator \cite{Koopman31}. Both methods use the snapshots of the fully discretized PDEs in time.  The POD solves a low dimensional model by Galerkin projection, whereas DMD is equation-free, the reduced solution is given in form of Fourier series in time and space. We would like to remark that due to its computational efficiency, the DMD is used for financial applications like in high-frequency trading \cite{Mann16} and stock market data analysis \cite{Cui16}.

The paper is organized as follows: In the next section, Section~\ref{HestonPDE}, we give the FOM of Heston's PDE applying dG discretization in space and backward Euler discretization in time. In Section~\ref{secrom}, the reduced-order modeling of Heston's PDE by POD and DMD is described. In the last section, Section~\ref{num}, we present numerical results for the butterfly spread, European call, and digital options by comparing the POD and DMD reduced solutions with respect to accuracy and speed-up. The paper ends with some conclusions.

\section{Full order model}
\label{HestonPDE}

\subsection{Heston model as diffusion-convection-reaction equation}

 The Heston model can be characterized by a two-dimensional diffusion-convection-reaction equation with variable coefficients.
  Let $u^S(t,v_t,S_t)$ be the value of a European option with the underlying price $S_t$ and volatilty $v_t$ and let $g(v_T,S_T)$ denote the payoff
  received at maturity $T$. Then, the option price $u^S(t,v_t,S_t)$ under Heston model satisfies the following linear two-dimensional variable coefficient diffusion-convection-reaction equation  ~\cite{Heston01041993}
\begin{equation}\label{pdeS1}
\frac{\partial u^S}{\partial t}+\mathcal{J}_{t}^{S}u^S-r_du^S=0,
\end{equation}
with the terminal condition
\begin{equation*}
u^S(T,v_T,S_T)=g(v_T,S_T),
\end{equation*}
where $v_t>0,$ $S_t>0,$ $t\in[0,T],$ and
\begin{equation*}
\mathcal{J}_{t}^{S} u^S=\frac{1}{2}S^2v\frac{\partial^2 u^S}{\partial S^2}+(r_d-r_f)S\frac{\partial u^S}{\partial S}
  +\rho\sigma Sv \frac{\partial^2 u^S}{\partial v\partial{S}}
  + \frac{1}{2}\sigma^2 v\frac{\partial^2 u^S}{\partial v^2}+\kappa(\theta-v)\frac{\partial u^S}{\partial v}.
\end{equation*}

Here, $r_d$ is the domestic interest rate, $r_f$ is the foreign interest rate, $\theta$ is the long-run mean level of $v_t$, $\sigma$ is the volatility of the volatilty and $\rho$ is the correlation coefficient. Applying the so called  log transformation $x=\log{(S/K)}$ and $\tau=T-t$ with $u(\tau,v,x)=u^S(T-\tau,v,Ke^x),$  PDE \eqref{pdeS1} is converted to the following equation
\begin{equation}\label{pdex1}
\frac{\partial u}{\partial \tau}- \mathcal{J}_{\tau}^x u+r_du=0,
\end{equation}
where $v>0,$ $x\in(-\infty,\infty),$ $\tau\in[0,T],$ and
\begin{align*}
\mathcal{J}_{\tau}^x u&=\frac{1}{2}v\frac{\partial^2 u}{\partial
  x^2}+(r_d-r_f-\frac{1}{2}v)\frac{\partial u}{\partial x}
  +\rho\sigma v \frac{\partial^2 u}{\partial
   v\partial{x}}
  +\frac{1}{2}\sigma^2 v\frac{\partial^2 u}{\partial
   v^2}+\kappa(\theta-v)\frac{\partial
    u}{\partial v}.\label{diffconvrecx}
\end{align*}
Note that due to the substitution $\tau=T-t,$ PDE (\ref{pdex1}) can also be regarded as a forward equation with
the following initial condition
\begin{equation*}
u^0:=u(0,v,x)=g(v,Ke^x).
\end{equation*}

We consider an open bounded domain $\Omega$ with the boundary $\Gamma = \Gamma_D\cup\Gamma_N$, where on $\Gamma_D$ the Dirichlet and on $\Gamma_N$  the Neumann boundary conditions are prescribed, respectively. Then, the log-transformed PDE given in (\ref{pdex1})  is expressed as a diffusion-convection-reaction equation
\begin{subequations}\label{convdiff_x}
\begin{align}
  \frac{\partial u}{\partial \tau}-\nabla\cdot (A \nabla u)+ b\cdot\nabla u+r_du&=0  & \text{in } (0,T]\times\Omega ,\label{convdiff_xmodel}\\
u(t,\bm{z}) &= u^D(t,\bm{z})  & \text{on } (0,T] \times\Gamma_D, \label{convdiff_xmodelDB}\\
A\nabla u(t,\bm{z})\cdot\bm{n} &= u^N(t,\bm{z})  & \text{on } (0,T]\times\Gamma_N,\label{convdiff_xmodelNB}\\
u(0,\bm{z}) &= u^0(\bm{z})  & \text{in } \{0\}\times\Omega ,\label{convdiff_xmodel_init}
\end{align}
\end{subequations}
where $\bm{n}$ is the outward unit normal vector, $\bm{z}=(v,x)^T$,  throughout this paper, denotes the spatial coordinates.
 In \eqref{convdiff_x}, the diffusion matrix and convective field  are given by
  \begin{eqnarray*}
  A  =
  \frac{1}{2}v\left( \begin{array}{ccc}
 \sigma^2  & \rho\sigma \\
  \rho\sigma  & 1
  \end{array} \right)\quad\text{and} \quad b =
  v \left( \begin{array}{c}
  \kappa\\
 \frac{1}{2}
\end{array}\right) +
  \left( \begin{array}{c}
  -\kappa\theta+\frac{1}{2}\sigma^2\\
 -(r_d-r_f)+ \frac{1}{2}\rho\sigma
\end{array}\right).
  \end{eqnarray*}
\begin{remark}
Although, the transformed PDE (\ref{pdex1}) is defined on the computational domain $(0,\infty)\times(-\infty,\infty)$, dGFEM  must be performed  on a bounded spatial region $\Omega=(v_{\text{min}},v_{\text{max}},)\times(x_{\text{min}},x_{\text{max}})$ for the numerical simulations, which is known as localization in option pricing models.
\end{remark}

\subsection{Variational form of Heston's model}
\label{variationalform}

We introduce the weak formulation of Heston's model as a parabolic  convection-diffusion-reaction equation \eqref{convdiff_xmodel}-\eqref{convdiff_xmodel_init}.
Let $L^2(\Omega)$ be the space consisting of all square integrable functions on $\Omega$, $H^1(\Omega)$ denote the Hilbert space of all functions having square integrable first-order partial derivatives, and $H_0^1(\Omega )=\{w\in H^1(\Omega): w=0\:\: \text{on} \:\:\Gamma_D\}$. The weak form of  \eqref{convdiff_x} is obtained by multiplying with a test function $w\in H_0^1(\Omega )$ and integrating by parts over the domain $\Omega$. Then, for a.e. $t\in(0,T]$, we seek a solution $u(t,v,x)\in H_D^1(\Omega )=\{w\in H^1(\Omega): w=u^D\:\: \text{on} \:\:\Gamma_D\}$ satisfying

\begin{subequations}\label{cont_weak}
\begin{align}
\int_{\Omega} \frac{\partial u}{\partial \tau}w d\bm{z} + a(u,w) &= \int_{\Gamma_N} u^Nw ds  &\forall w\in H_0^1(\Omega ),\\
\int_{\Omega} u(0,\bm{z})wd\bm{z} &=  \int_{\Omega}u^0 wd\bm{z}  &\forall w\in H_0^1(\Omega ),
\end{align}
\end{subequations}
where $ds$ is the arc-length element on the boundary. In  \eqref{cont_weak}, $a(u,w)$ is the bilinear form given by
$$
a(u,w)=\int_{\Omega}\left( A \nabla u\cdot \nabla w + b \cdot\nabla uw+r_duw\right)d\bm{z}, \quad \forall w \in H_0^1(\Omega).
$$

We assume that the matrix $A$ is positive definite for $v >0$ and $\rho \in (-1,1)$ which is usually satisfied.

\subsection{Discontinuous Galerkin discretization in space}
\label{sec:dg}

The symmetric interior penalty Galerkin (SIPG) method is the commonly used dG method, which
enforces boundary conditions weakly \cite{riviere08dgm}. Let the mesh $\xi_{h} =\{K\}$ be a partition of the domain $\Omega$ into a family of shape regular elements (triangles). We set the mesh-dependent finite-dimensional solution and test function space by
$$
W_{h}=W_{h}(\xi_{h})=\left\{ w\in L^2(\Omega ) : w|_{K}\in\mathbb{P}_k(K) ,\; \forall K\in \xi_{h} \right\}\not\subset H_0^1(\Omega ),
$$
where the functions in $W_{h}$ are discontinuous along the inter-element boundaries. These discontinuities lead to the fact that on an interior edge $e$ shared by two neighboring triangles $K_i$ and $K_j$, there are two different traces from either triangle. Thus, for convenience, we define the jump and average operators of a function $w\in W_{h}$ on $e$ by
$$
\jump{w}:= w|_{K_i}\bm{n}_{K_i} + w|_{K_j}\bm{n}_{K_j}\; , \quad \av{w}:=\frac{1}{2}(w|_{K_i} + w|_{K_j}).
$$
On a boundary edge $e\subset\partial\Omega$, we set $\jump{w}:= w|_{K}\bm{n}$ and $\av{w}:=w|_{K}$. In addition, we form the sets of inflow and outflow edges by
$$
\Gamma^- = \{ \bm{z}\in\partial\Omega : \bm{b}(v)\cdot\bm{n}(v,x)<0\} \; , \qquad \Gamma^+ = \partial\Omega\setminus\Gamma^- ,
$$
$$
\partial K^- = \{ \bm{z}\in\partial K: \bm{b}(v)\cdot\bm{n}_K(v,x)<0\} \; , \qquad \partial K^+ = \partial K\setminus\partial K^-,
$$
where $\bm{n}_K$ denotes the outward unit vector on an element boundary $\partial K$. Moreover, we denote by $\Gamma_{h}^0$ and $\Gamma_{h}^D$ the sets of interior and Dirichlet boundary edges, respectively, so that the union set is $\Gamma_{h}=\Gamma_{h}^0\cup\Gamma_{h}^D$. Then, in space SIPG discretized semi-discrete system of the PDE \eqref{convdiff_x} reads as: for a.e. $t\in (0,T]$, for all $w_h\in W_{h}$, find $u_h:=u_h(t,\bm{z})\in W_{h}$ such that
\begin{subequations} \label{dg}
\begin{align}
\int_{\Omega}\frac{\partial u_{h}}{\partial t}w_{h}d\bm{z} + a_{h}(t;u_{h},w_{h}) &=l_{h}(w_{h}) ,\\
\int_{\Omega} u_{h}(0,\bm{z})w_hd\bm{z} &= \int_{\Omega} u^0w_hd\bm{z},
\end{align}
\end{subequations}
with the (bi)linear forms:
\begin{align*}
a_{h}(t;u_{h}, w_{h})=& \sum \limits_{K \in {\xi_{h}}} \int_{K} \left( A \nabla u_{h}\cdot\nabla w_{h} + b \cdot \nabla u_{h} w_{h} + r_du_hw_h\right)d\bm{z} \\
& + \sum \limits_{ e \in \Gamma_{h}} \int_{e} \left( \frac{\sigma_e}{h_{e}}\jump{u_{h}}\cdot\jump{w_{h}} - \av{A \nabla w_{h}}\jump{u_{h}} - \av{A  \nabla u_{h}} \jump{w_{h}}\right)ds  \\
& + \sum \limits_{K \in {\xi_{h}}} \left( \int\limits_{\partial K^-\setminus\partial\Omega } b \cdot \bm{n}_K (u_{h}^{out}-u_{h})  w_{h} ds -  \int\limits_{\partial K^-\cap \Gamma^{-}} b\cdot \bm{n}_K u_{h} w_{h} ds \right), \\
l_{h}( w_{h})=&  \sum \limits_{e \in {\Gamma_{h}^N}} \int_{e} u^N w_{h} d\bm{z} + \sum \limits_{e \in {\Gamma_{h}^D}} \int_{e} u^D\left( \frac{\sigma_e}{h_{e}}w_{h} - A \nabla w_{h} \right)ds\\
&  - \sum \limits_{K \in {\xi_{h}}}\int\limits_{\partial K^-\cap \Gamma^{-}} b\cdot \bm{n}_K u^D w_{h} ds,
\end{align*}
where $u_{h}^{out}$ denotes the trace of $u_h$ on an edge $e$ from outside the triangle $K$. Here, the parameter $\sigma_e$ is called the penalty parameter, and it should be selected sufficiently large to ensure the coercivity of the bilinear form \cite[Sec. 27.1]{riviere08dgm}.

The solution of  SIPG semi-discretized Heston model \eqref{dg} is given as
\begin{equation}\label{4}
u_h(t,\bm{z})=\sum^{n_e}_{m=1}\sum^{n_k}_{j=1}u^{m}_{j}(t) \varphi^{m}_{j}(\bm{z}),
\end{equation}
where $\varphi^{m}_{j}$ and $u^{m}_{j}$, $j=1, \ldots, n_k$, $m=1, \ldots, n_e$, are the basis functions spanning the space $W_{h}$ and the unknown coefficients, respectively. The number $n_k$ denotes the local dimension of each dG element with  $n_k=(k+1)(k+2)/2$ for 2D problems, and $n_e$ is the number of dG elements (triangles). Substituting \eqref{4} into \eqref{dg} and choosing $\upsilon=\varphi^{k}_{i}, \: i=1, \ldots, n_k$, $k=1, \ldots, n_e$, we obtain the following semi-linear system of ordinary differential equations (ODEs)
for the unknown coefficient vector $\bm{u}:= \bm{u}(t) =(u_{1}^{1}(t), \ldots , u_{n_k}^{1}(t),\ldots , u_{1}^{n_e}(t), \ldots , u_{n_k}^{n_e}(t))^T\in\mathbb{R}^{{\mathcal N}}$
\begin{equation}\label{fom}
\bm{M}\bm{u}_{t} + \bm{A}\bm{u}   = \bm{l},
\end{equation}
where $\bm{M}$ is the mass matrix and  $\bm{A}$ is the stiffness matrix, with the entries $(\bm{M})_{ij}=(\varphi^{j}, \varphi^{i})_{\Omega}$ and $(\bm{A})_{ij}=a_h(\cdot ; \varphi^{j}, \varphi^{i})$, $1\leq i,j \leq {\mathcal N}:=n_k\times n_e$, and $\bm{l}$ is the linear right-hand side vector with entries $(\bm{l})_{i}=l_h(\varphi^{i})$. For the time discretization, we consider a subdivision of $[0,T]$ into $J$ time intervals $I_n=(t_{n-1},t_{n}]$ of length $\Delta t$, $n=1,2,\ldots , J$, with $t_0=0$. Then, the backward Euler solution of the fully discrete formulation of \eqref{dg} reads as: for $t=0$ set $u_h^0\in W_{h}$ as the projection (orthogonal $L^2$-projection) of $u^0$ onto $W_{h}$ and for $n=1,2,\ldots , J$, find $u_h^n:=u_h^n(\bm{z})\in W_{h}$ satisfying for all $w_h\in W_{h}$
\begin{equation*}
\int_{\Omega}\frac{u_{h}^{n+1}-u_{h}^{n}}{\Delta t}w_{h}d\bm{z} + a_{h}(t_{n+1};u_{h}^{n+1},w_{h}) =l_{h}(w_{h}),
\end{equation*}
or in matrix-vector form
\begin{equation}\label{fullydisc}
(\bm{M} + \Delta t \bm{A} )\bm{u}^{n+1}   =  \bm{M}\bm{u}^n + \Delta t\bm{l}^{n+1}.
\end{equation}
 The coefficient matrix $(\bm{M} + \Delta t \bm{A} )$ is factorized by LU decomposition at the initial time step and used in all successive time steps.

\section{Reduced-order modeling}
\label{secrom}

Both POD and DMD are snapshot-based post-processing
algorithms which may be applied equally well to data obtained in simulations or in experiments.
The POD is based on the observation that the dynamics of the PDE are optimally contained in a small number of modes computed from a singular value decomposition \cite{Kunisch01}.  After determination of the POD basis by a predetermined cut-off value, the truncated POD modes are used as the basis for Galerkin expansion for the reduced-order dynamical system.
The DMD can be interpreted as a model order reduction technique like the POD with temporal POD modes \cite{Tu14}.
While POD modes are characterized by spatial orthogonality and with multi-frequential temporal content, DMD modes are non-orthogonal but each of them possesses a single temporal frequency. 
Both methods are conceptually different; POD attempts to build a low-dimensional basis for the
solution and DMD attempts to build a low-dimensional basis for the Koopman
operator \cite{Koopman31}.

In the DMD,
the eigenfunctions of an unknown linear time-independent operator are approximated, which
can be thought of as a finite-dimensional approximation of the infinite-dimensional
Koopman operator \cite{Koopman31}. The dimensionality reduction occurs through the  approximation of the
infinite-dimensional set of Koopman modes with a finite-dimensional set of
eigenvectors \cite{Tu14}.
Because there are time dynamics associated with each
eigenfunction, no analog of the Galerkin projection is needed to use the DMD
modes. DMD is equation-free, where the solutions are given in form of Fourier series in space and time. This feature of DMD allows making future predictions.
However, the DMD modes are not orthogonal and it requires in general more modes than
an equally accurate approximation of the data with POD.

\subsection{Proper orthogonal decomposition}

To form a ROM via POD, we construct a low-dimensional ODE system by a Galerkin projection procedure. The FOM \eqref{fom} is of a large dimension ${\mathcal N}$. The ROM of a small dimension $N\ll {\mathcal N}$ is of the form:
\begin{equation}\label{rom}
\partial_t\bm{u}_{r} + \bm{A}_r\bm{u}_{r} = \bm{l}_r,
\end{equation}
with the reduced stiffness matrix $\bm{A}_r=\bm{U}^T\bm{A}\bm{U}$, the reduced  vector $\bm{l}_r=\bm{U}^T\bm{l}$, and the unknown vector $\bm{u}_{r}\in\mathbb{R}^N$ is the coefficient vector of the reduced solution $u_{h,r}$ from the $N$ dimensional reduced space $W_{h,r}$ spanned by the reduced basis functions $\{U_i\}_{i=1}^N$. The matrix $\bm{U}=[\bm{U}_{1},\ldots ,\bm{U}_{N}]\in\mathbb{R}^{{\mathcal N}\times N}$ is the matrix whose columns $\bm{U}_{i}$'s will be called as the reduced basis modes, and they are the coefficient vectors of the reduced basis functions $U_{i}$'s. Indeed, we have the following relations
\begin{equation*}
u_{h,r}(t,\bm{z})= \sum_{i=1}^{N} (\bm{u}_{r})_i(t)U_{i}(\bm{z}), \quad U_{i}(\bm{z}) = \sum_{j=1}^{{\mathcal{N}}} \bm{U}_{j,i}\varphi_{j}(\bm{z}),\quad \bm{u}= \bm{U} \bm{u}_r,
\end{equation*}
where $\bm{u}$ is the coefficient vector of the FOM solution, and $\{\varphi_{i}\}_{i=1}^{\mathcal{N}}$ are the dG finite elements basis functions, for which then the reduced basis functions $U_i$'s lies in the space $W_h$ and so $W_{h,r}\subset W_h$. The ROM \eqref{rom} is obtained by the projection of the FOM \eqref{fom} by the matrix $\bm{U}$ and by the substitution $\bm{u}= \bm{U}\bm{u}_r$. Another fact with \eqref{rom} is the M-orthogonality of the reduced modes $\bm{U}_i$'s, i.e. $\bm{U}_i^TM\bm{U}_j=\delta_{ij}$ where $\delta_{ij}$ is the Kronecker Delta. For this reason, the reduced mass matrix in \eqref{rom} is the identity matrix ($\bm{U}^TM\bm{U} =\bm{I}$).

The computation of the reduced modes $\bm{U}_i$'s are based on the fact that the reduced basis functions $U_i$'s are the solution of the optimization problem \cite{Volkweinn13}
\begin{equation}\label{optim}
\begin{aligned}
\min_{U_1,\ldots ,U_N} \frac{1}{J}\sum_{j=1}^J \left\| u_h^j - \sum_{i=1}^N \overbrace{(u_h^j,U_i)_{L^2(\Omega)}}^{(\bm{u}_{r})_i}U_i\right\|_{L^2(\Omega)}^2 \\
\text{subject to } (U_i,U_j)_{L^2(\Omega)} = \bm{U}_i^TM\bm{U}_j=\delta_{ij} \; , \; 1\leq i,j\leq N,
\end{aligned}
\end{equation}
where $u_h^j:=u_h^j(\bm{z})\approx u_h(t_j,\bm{z})$ is the solution through the FOM at $t=t_j$. The minimization problem \eqref{optim},  on the other hand, is equivalent to the eigenvalue problem
\begin{equation*}\label{eigenprob}
\widehat{\mathcal{S}}\widehat{\mathcal{S}}^T\widehat{\bm{U}}_{\cdot ,i}=\sigma_{i}^2\widehat{\bm{U}}_{\cdot ,i}, \quad i=1,2,\ldots ,N
\end{equation*}
where $\widehat{\bm{U}}_{\cdot, i}=R\bm{U}_{\cdot ,i},$ $\widehat{\mathcal{S}}=R\mathcal{S}$, $R^T$ is the Cholesky factor of the mass matrix $M$, and $\mathcal{S}=[\bm{u}^1,\ldots ,\bm{u}^J]\in\mathbb{R}^{\mathcal{N}\times J}$ is the snapshot matrix, where $\bm{u^i}\approx\bm{u}(t_i),$ $i=1,2,\ldots ,J$. To have an enough accuracy, the number $N$ of POD modes to be used is selected as the first integer satisfying the relative information content $I(N)$ criteria:
\begin{equation*}
I(N)=\frac{\sum_{i=1}^{N}\sigma_i^2}{\sum_{i=1}^{s}\sigma_i^2}\geq 1-\varepsilon^2,
\end{equation*}
where $I(N)$ represents the energy captured by the first $N$ POD modes for a given tolerance $\varepsilon$, and $s$ is the column rank of the snapshot matrix $\mathcal{S}$.

\subsection{Dynamic mode decomposition}

We consider the dynamic system  on a manifold $\mathcal{M}$
\[{\bm u}^{n+1}=g({\bm u}^n).\]
The linear infinite dimensional Koopman operator \cite{Koopman31} maps any function $f:\mathcal{M}\rightarrow \mathbb{C}$ into
 \begin{equation*}
 {\mathcal  A } f({\bm u}^n)=f(g({\bm u}^n)).
  \end{equation*}
Then, the function $f({\bm u}^n)$ can be  written as
$$f({\bm u}^n) = \sum_{j=1}^{\infty} \alpha_j v_j \lambda_j^n,$$
where the eigenvectors $v_j\in \mathbb{C}$ denote the dynamic modes as Koopman eigenfunctions,  $\lambda_j\in \mathbb{C}$ denote the Ritz eigenvalues of the eigenvalue problem ${\mathcal A}v_j  = \lambda v_j$ and $\alpha_j\in \mathbb{C}$ denotes the amplitudes  of the Koopman modes.

The DMD represents the eigendecomposition \cite{Tu14} of an approximating linear operator $\mathcal A$ corresponding to the Schmidt operator \cite{Schmid10}, which a is a special case of
Koopman operator acting on the dynamic variable ${\bm u}$.

We consider the snapshot matrix $\mathcal{S}=[\bm{u}^1,\ldots ,\bm{u}^J]$ in $\mathbb{R}^{\mathcal{N} \times J}$ with $$\mathcal{S}_0=[\bm{u}^1,\ldots ,\bm{u}^{J-1}], \quad  \mathcal{S}_1=[\bm{u}^2,\ldots ,\bm{u}^J],$$
as time discrete solutions of \eqref{fullydisc}.

 The DMD is based on the fact that for sufficiently large set of snapshots, the snapshots can be written as
 \begin{equation*}
 {\bm u}^{n+1}=\mathcal {A} {\bm u}^n
 \end{equation*}
  where $\mathcal {A}$ is the Koopman operator. Hence, the snapshots form a Krylov sequence ${\mathcal S}_1 = \{ {\bm u}^1, {\mathcal A}{\bm u}^1, {\mathcal A}^2{\bm u}^1, \cdots, {\mathcal A}^{J-1}{\bm u}^1 \}$, where the following snapshots became linearly dependent on the previous ones. The last snapshot ${\bm u}^J$ is expressed with the error term  $R$ as:
$$
{\bm u}^J = c_1{\bm u}^1 + c_2{\bm u}^2+ \cdots + c_{J-1}{\bm u}^{J-1} + R
$$

The DMD algorithms are based on the minimization of the norm $R$. This can be expressed  equivalently as
$
\mathcal{S}_1\approx {\mathcal A }\mathcal{S}_0,
$
so that $\tilde{\mathcal A}$ is the minimizer of
\begin{equation*}
\|\mathcal{S}_1- \tilde{\mathcal A} \mathcal{S}_0\|_F
\end{equation*}
in the Frobenious norm.

Different DMD algorithms are developed for the estimation of Koopman modes, eigenvalues and amplitudes from the given set of snapshots. In this paper we consider the  exact DMD algorihm in \cite{Tu14} and a variant of DMD algorithm in \cite{Chen12}.

 \begin{algorithm}
 \caption{Exact DMD Algorithm (Tu et al.~\cite{Tu14} ) }\label{exactDMDalg}
\begin{flushleft}
 \textbf{Input:} Snapshots $\mathcal{S}=[\bm{u}^1,\ldots ,\bm{u}^J]$ with $\bm{u^i}\approx\bm{u}(t_i).$
  \\
 \textbf{Output:} DMD modes $\Phi^{\text{DMD}}.$
\end{flushleft}
 \begin{algorithmic}[1]
 \STATE Define the matrices $\mathcal{S}_0=[\bm{u}^1,\ldots ,\bm{u}^{J-1}]$ and $\mathcal{S}_1=[\bm{u}^2,\ldots ,\bm{u}^J].$
 \STATE Take SVD of $\mathcal{S}_0:$ $\mathcal{S}_0= {\bm U}{\bm \Sigma}{\bm V}^*.$
 \STATE Consider $\tilde{\mathcal A}={\bm U}^*\mathcal{S}_1 {\Sigma}^{-1}{\bm V}.$
 \STATE Determine eigenvalues and eigenvectors of $\tilde{\mathcal A}W=\Lambda W$
 \STATE Obtain $\Phi^{\text{DMD}}=\mathcal{S}_1{\bm \Sigma}^{-1}{\bm V} W.$
  \end{algorithmic}
 \end{algorithm}

	\begin{algorithm}
		\caption{Variant of DMD Algorithm (Chen et al.~\cite{Chen12})  }\label{exactDMDalg_Chen}
		\begin{flushleft}
		\textbf{Input:} Snapshots $\mathcal{S}=[\bm{u}^1,\ldots ,\bm{u}^J]$ with $\bm{u^i}\approx\bm{u}(t_i).$
		\\
		\textbf{Output:} DMD modes $\Phi^{\text{DMD}}.$
		\end{flushleft}
		\begin{algorithmic}[1]
			\STATE Define the matrices $\mathcal{S}_0=[\bm{u}^1,\ldots ,\bm{u}^{J-1}]$ and $\mathcal{S}_1=[\bm{u}^2,\ldots ,\bm{u}^J].$
			\STATE Take SVD of $\mathcal{S}_0^T\mathcal{S}_0:$ $\mathcal{S}_0^T\mathcal{S}_0= {\bm V} {\bm \Sigma}^2{\bm V}^*.$
			\STATE Let ${\bm U} =\mathcal{S}_0{\bm V}{\bm \Sigma}^{-1}.$
			\STATE Consider $\tilde{\mathcal A}={\bm U}^*\mathcal{S}^1 {\bm \Sigma}^{-1}V$
			\STATE Determine eigenvalues and eigenvectors of $\tilde{\mathcal A}W=\Lambda W$
			\STATE Obtain $\Phi^{\text{DMD}}=\mathcal{S}^1{\bm \Sigma}^{-1}{\bm V} W.$
		\end{algorithmic}
	\end{algorithm}

In POD the modes are ranked by energy level through the POD singular values. There is no such criteria for ranking the contributions of the different DMD modes. Different criteria are developed depending on what can be considered important for the models used. The DMD modes can then be selected based on their amplitude or based on their frequency/growth rate. The amplitude criterion is also not perfect because there exist modes with very high amplitudes but which are very fast damped. The selection based on frequency/growth rate has also disadvantages because it relies on a priori physical knowledge.
Additionally, spatial non-orthogonality of the DMD modes
may introduce a poor quality of approximation
when only a subset of modes with the largest amplitude is retained. Recently several algorithms are developed for selecting optimal amplitudes and extracting the
desired frequencies, spatial profiles using combinatorial search \cite{Chen12}.

We give here  briefly the optimal selection of amplitudes of
extracted DMD modes following \cite{Javanovich14}. The dynamics of the reduced system in $r$ dimensional subspace is governed by
$$
{\bm u}^{n+1} = \tilde {\mathcal A } {\bm u}^n,
$$
where $\tilde{\mathcal A }= {\bm U}^*\mathcal{S}_1{\bm V} {\bm \Sigma}^{-1}$ is the reduced matrix obtained from the DMD algorithms Algorithm~\ref{exactDMDalg} and Algorithm~\ref{exactDMDalg_Chen}. The reduced solution can be written as linear combination of DMD modes
$$
{\bm u}_r = \sum_{i=1}^r \alpha_i\phi_i\lambda_i^k, \quad k\in [0,1,\ldots,J-1]
$$
The unknown optimal amplitudes ${\bm \alpha} = (\alpha_1, \cdots, \alpha_r)$ are then determined by solving following minimization problem \cite{Javanovich14}

 \begin{equation*}
 \min_{\substack{\alpha}}\|\mathcal{S}-\Phi^{\text{DMD}}D_{\alpha}V_{\text{and}}\|_{F}^2,
 \end{equation*}
 where $D_{\alpha}=\text{diag}(\alpha_1,\ldots,\alpha_r)$ and $V_{\text{and}}$ is the Vandermonde matrix $V(\gamma_1,\ldots,\gamma_r)$ with $\gamma_i=\exp(\omega_ i t),\; i=1,2,\ldots ,k.$
Let
$P=(W^* W)\circ(\overline{V_{\text{and}} V_{\text{and}}^*}),\; q=\overline{\text{diag}(V_{\text{and}} V\Sigma^* W)}$
where $^*$ denotes the conjugate transpose, $\circ$ denotes the elementwise multiplication.
 Then, ${\bm \alpha}_{\text{opt}}=P^{-1}q$.

Finally the approximate DMD solution is given as
\begin{equation*}
{\bf u}_r(t)=\sum_{j=1}^{r}\alpha_j(0)\phi_j(z) \exp(\omega_ j t) = \Phi^{\text{DMD}}\text{diag}(\exp(\omega t)) {\bm \alpha}_{Opt}(0),
\end{equation*}
where
$$\Phi^{\text{DMD}}=[\phi_1,\ldots,\phi_k ],\quad {\bm  \alpha}_{Opt}(0)=[\alpha_1(0),\ldots, \alpha_k(0)], \quad \omega_j=\log{(\lambda_j)}/\Delta t.
$$
The initial amplitude of the modes are determined by ${\bm \alpha}_{Opt}(0) = \Phi^\dag {\bm u}^1 $, where
$\Phi^\dag $  denotes the Moore-Penrose pseudo inverse of the DMD modes $\Phi^{\text{DMD}}.$

\section{Numerical results}
\label{num}

In this section, we present numerical results for different options to compare the reduced solutions of POD and DMD with respect to accuracy and speed-up.
 In all numerical tests, we have used linear dG elements in space and backward Euler in time. For the computation of the DMD modes, we used the MATLAB Toolbox Koopman mode decomposition \cite{Budisic15}. The numerical simulations given in this paper are performed on
Windows 7 with an Intel Core i7, 2.9Ghz, and 8GB RAM using MATLAB R2014a.

\subsection{European  call option}

We impose the initial and  boundary conditions in~\cite{Winkler01}:
\begin{eqnarray}
u(\tau,v_{\text{min}},x)&=&Ke^{x-r_f\tau}\Phi(d_{+})-Ke^{r_d\tau}\Phi(d_{-})\nonumber\\
u(\tau,v_{\text{max}},x)&=&Ke^{x-r_f\tau}\nonumber\\
u(\tau,v,x_{\text{min}})&=&\lambda u(\tau,v_{\text{max}},x_{\text{min}})+(1-\lambda)u(\tau,v_{\text{min}},x_{\text{min}})\nonumber\\
\frac{\partial}{\partial \nu}u(\tau,v,x_{\text{max}})&=&A\nabla u\cdot\vec{\bm n}=\frac{1}{2}v K e^{x-r_f\tau}    \nonumber\\
u(0,v,x)&=&(Ke^x-K)^{+}\nonumber
\end{eqnarray}
where $\vec{\bm n}$ is the outward normal vector,
\begin{equation*}
d_{+}=\frac{x+\left(r_d-r_f+\frac{1}{2}v_{\text{min}}\right)\tau}{\sqrt{v_{\text{min}}\tau}},\quad d_{-}=\frac{x+\left(r_d-r_f-\frac{1}{2}v_{\text{max}}\right)\tau}{\sqrt{v_{\text{max}}\tau}}
\end{equation*}
and $\Phi(x)$ is the cumulative distribution function given as
\begin{equation*}
\Phi(x)=\frac{1}{2\pi}\int_{-\infty}^{x}e^{-y^2/2}dy.
\end{equation*}
We take the discretization parameters as $\Delta t=0.01,$ $N_x=96,$ and $N_v=48$  in the domain $[0.0025,0.5]\times[-5,5]$.
The parameter set for the European call option is taken from~\cite{Burkovska15} with strong negative correlation  $\rho=-0.9.$

\begin{table}[h]%
		\centering
		\caption{Parameter set for the European call option.}%
		\begin{tabular*}{425pt}{@{\extracolsep\fill}cccccccccc@{\extracolsep\fill}}
			\hline
			$\bm{\kappa}$ & $\bm{\theta}$  & $\bm{\sigma}$  & $\bm{\rho}$  & $\bm{r_d}$  & $\bm{r_f}$  & $\bm{T}$  & $\bm{S_0}$  & $\bm{K}$  & $\bm{v_0}$ \\
			\hline
			2.5 & 0.06 & 0.4 &  -0.9 & 0.0198 & 0 &  1& 1& 1&0.1683  \\
			\hline
		\end{tabular*}
	\end{table}

\begin{figure}[htb!]
\centering
\includegraphics[width=.45\columnwidth]{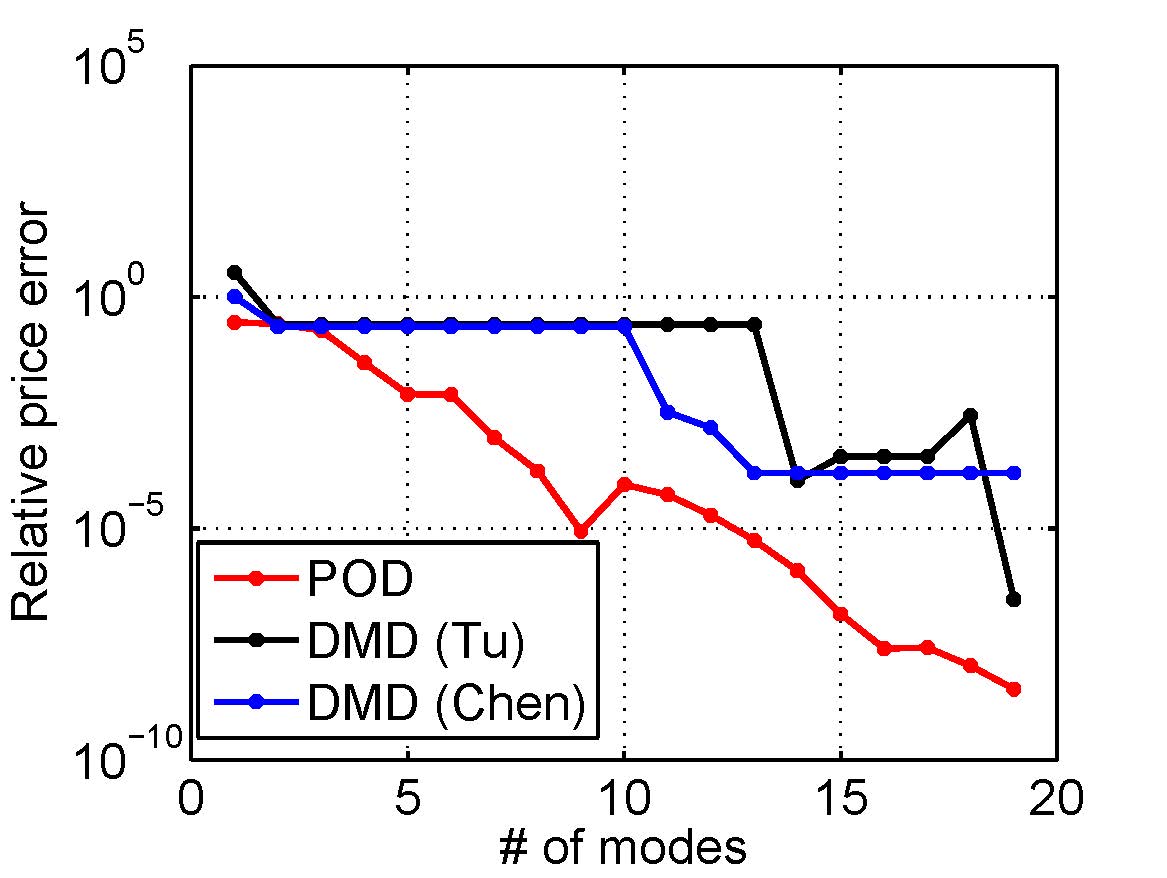}
\includegraphics[width=.45\columnwidth]{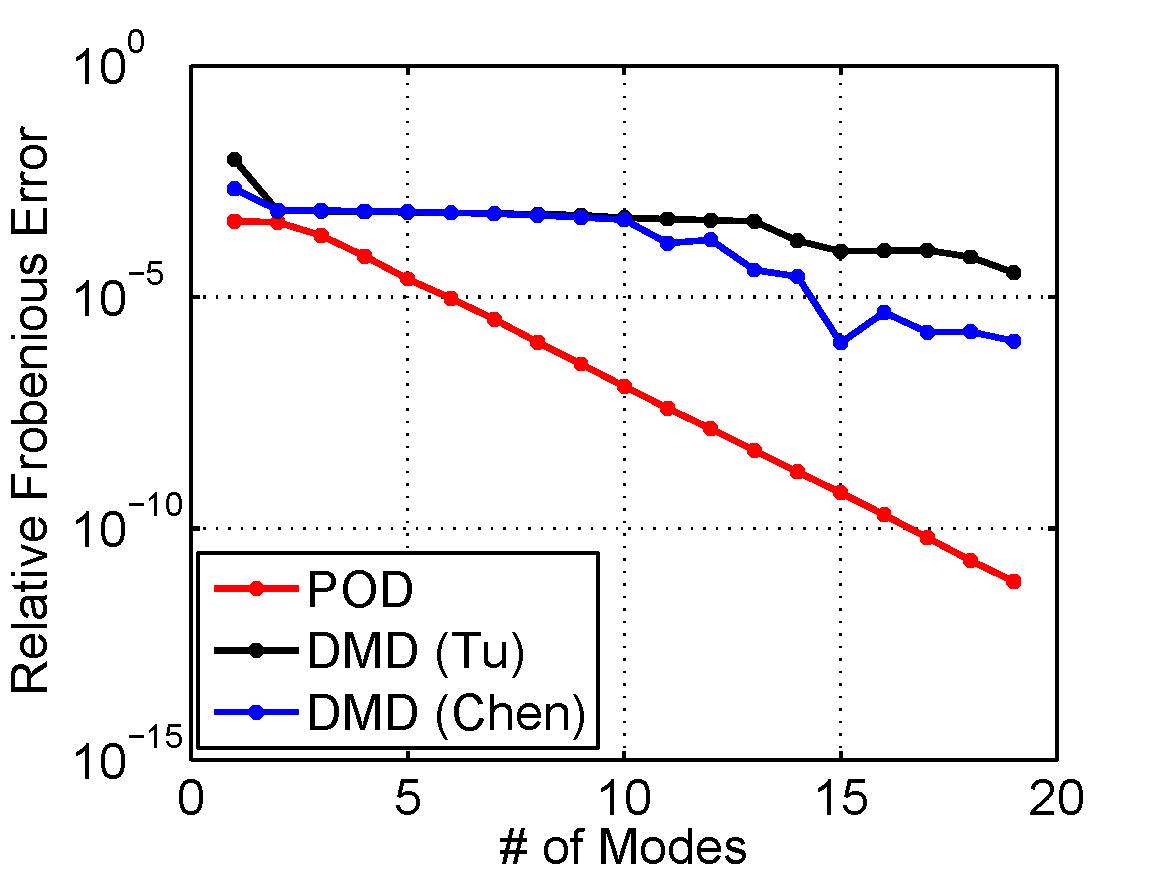}
\caption{Relative price error for $v_0=0.1683$ and $S_0=1$ (left), relative Frobenious error  (right)\label{fig:european1}}
\end{figure}

The relative price and Frobenious errors between the FOM and ROM solutions decay monotonically in \figurename~\ref{fig:european1} for the POD modes, whereas they reach plateaus and decrease more slowly for the DMD modes.

\begin{figure}[htb!]
\centering
\subfloat{\includegraphics[width=.35\columnwidth]{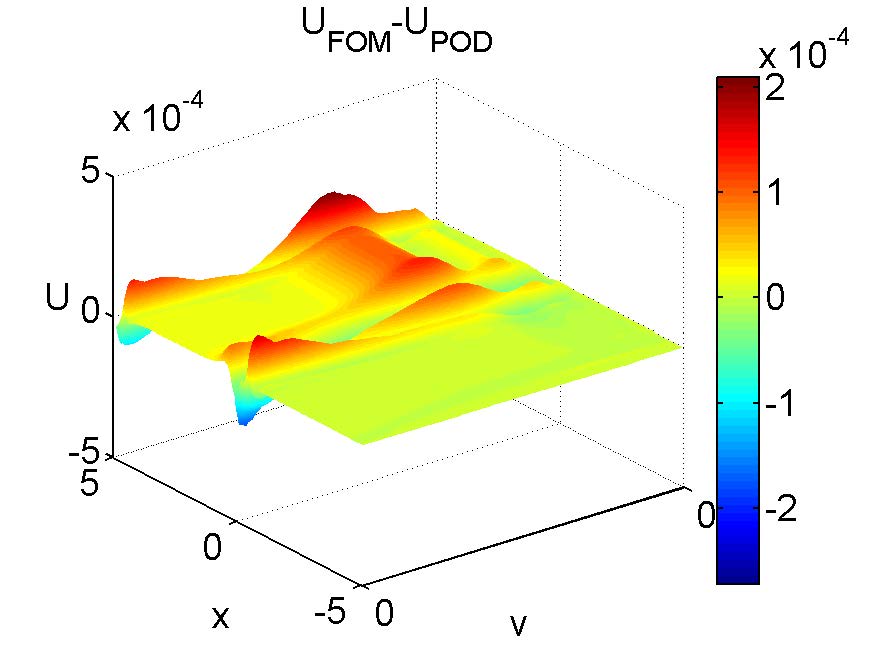}}
\subfloat{\includegraphics[width=.35\columnwidth]{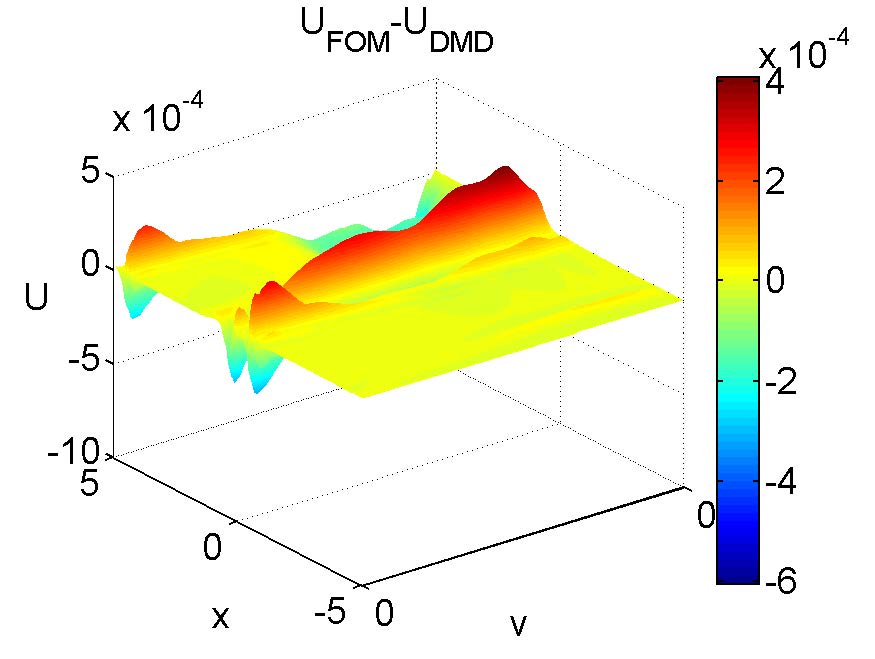}}
\subfloat{\includegraphics[width=.35\columnwidth]{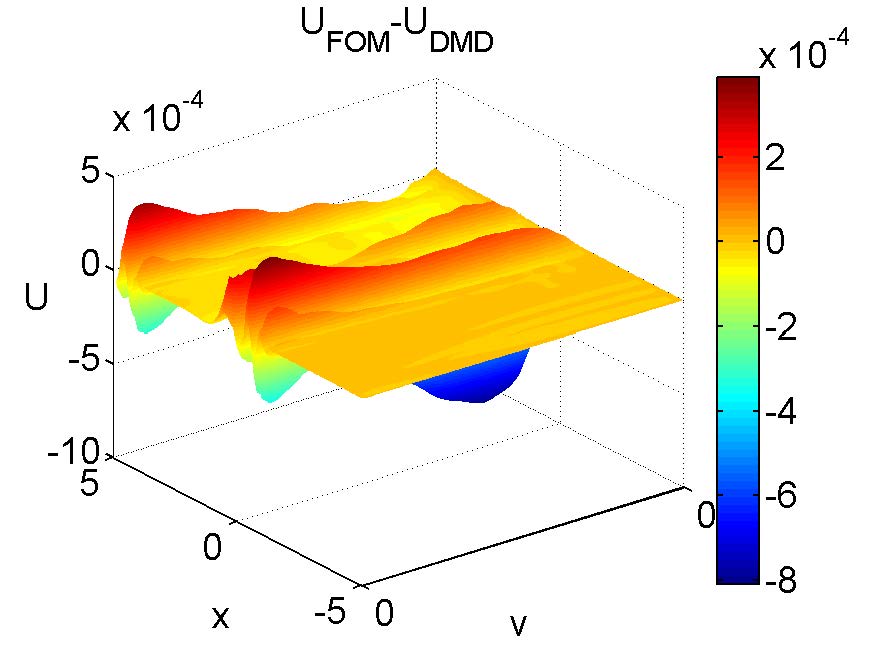}}
\caption{ROM-FOM errors: 8 POD modes (left), 12 DMD (Chen) modes (middle), 18 DMD (Tu) modes function (right).}
\label{fig:european2}
\end{figure}

In \figurename~\ref{fig:european2} the ROM-FOM errors are plotted at the almost same accuracy level for different POD and DMD modes. As expected the POD requires fewer modes than both DMD algorithms. As can be seen, in the neighborhood of  $x=0$ (i.e $S=K$) and at the boundaries relatively large errors are observed. One of the main reasons for this is the initial function $(Ke^x-K)^{+}$ has a discontinuity in its first derivative at $x=0.$
Hence, we may conclude that the reduced-order models could not resolve the full order European call option pricing problem for the at-the-money options.

\begin{figure}[htb!]
\centering
\includegraphics[width=.5\columnwidth]{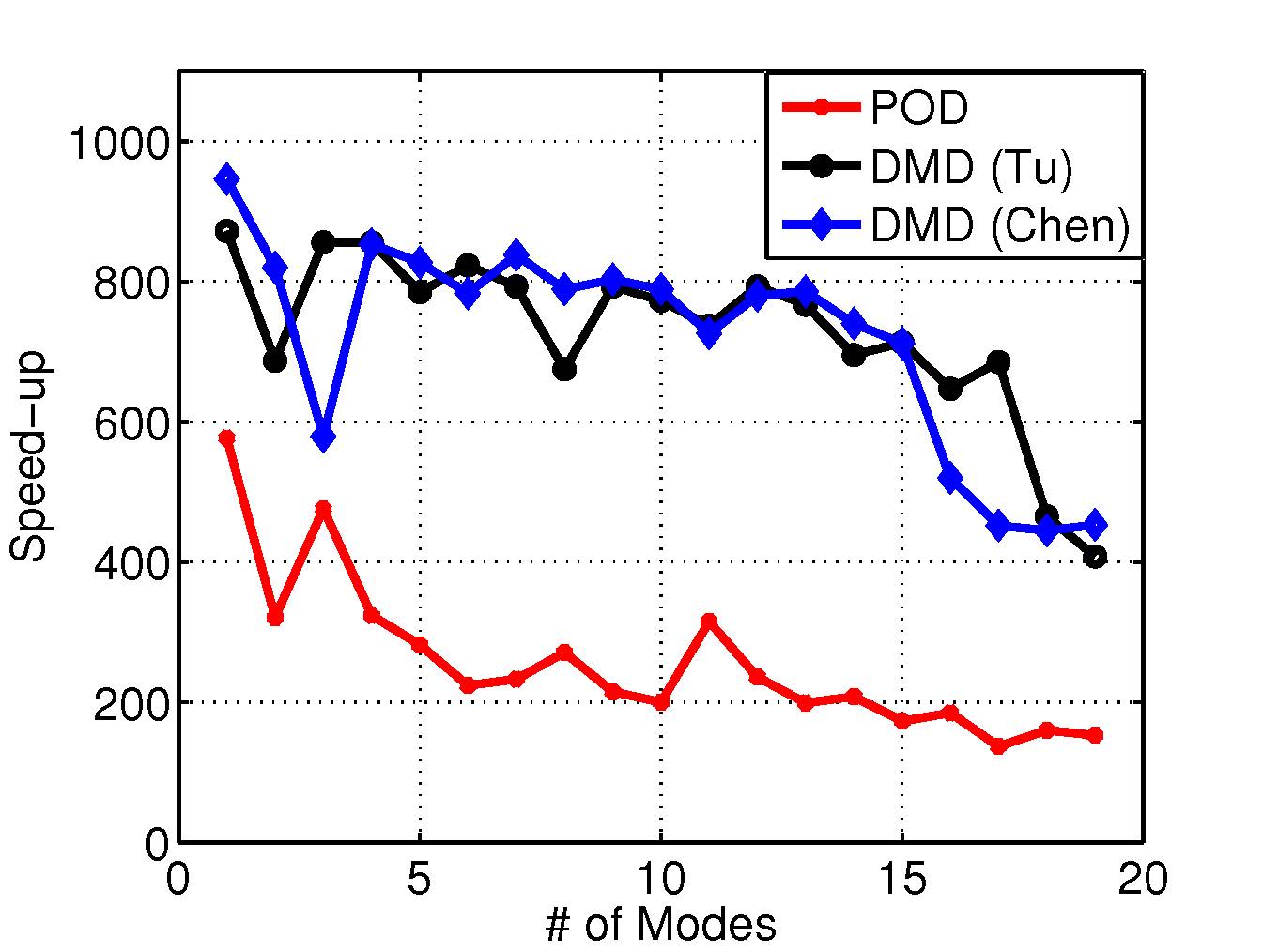}
\caption{Speed-up factors}
\label{fig:speed}
\end{figure}

The performance of the DMD over POD is clearly seen in Figure \ref{fig:speed}. The speed-up factors of both DMD algorithms increase more rapidly than of the POD with an increasing number of modes. This is due to the fact that the DMD produces equation-free solutions, whereas for the POD, for an increasing number of modes larger reduced-order Galerkin projected ordinary differential equations have to be solved.

\subsection{Butterfly spread}

The butterfly spread is composed of three call options with different strike prices.  Two call options are bought for $K_1<K_3,$ with a strike price $K_1$ and $K_3,$ and  two call options sold with a strike price $K_2=(K_1+K_3)/2.$  The payoff function is then given by
\begin{equation*}
    g(v,Ke^x)=(K_2 e^x-K_1)^+ - 2(K_2 e^x-K_2)^+ + (K_2 e^x-K_3)^+,
\end{equation*}
with  $x=\log{(S/K_2)}.$ Note that the payoff function (initial data) is non-differentiable at the strike prices $K_1$, $K_2$ and $K_3.$

Let $u(\tau,v,x)$ be the price of a butterfly spread option satisfying Heston's PDE (\ref{convdiff_x}) with $x=\log{(S/K_2)}$. We impose the homogeneous Dirichlet boundary conditions in the $x$-direction and homogeneous Neumann boundary conditions in the $v$-direction \cite{Ballestra13}:
$$
u(\tau,v,x_{\text{min}})=0, \quad u(\tau,v,x_{\text{max}})=0,
$$
$$
 \frac{\partial}{\partial v}u(\tau,v_{\text{min}},x) = 0\; , \quad
 \frac{\partial}{\partial v}u(\tau,v_{\text{max}},x) = 0,
$$
with an initial condition
    $$ u(0,v,x)= (K_2 e^x-K_1)^+ - 2(K_2 e^x-K_2)^+ + (K_2 e^x-K_3)^+.$$

The computational domain is taken as  for the European call option $(0.0025,0.5)\times(-5,5).$
The discretization parameters are taken as $\Delta t=0.01,$ $N_x=96,$ and $N_v=48$  and $K=0.5$, $K_1=0.1$ and $K_2=0.9$.
 The parameter set is the same as in ~\cite{Burkovska15} with positive correlation  $\rho=0.55$

\begin{table}[htb!]%
		\centering
		\caption{Parameter set for the butterfly spread.}%
		\begin{tabular*}{425pt}{@{\extracolsep\fill}ccccccccc@{\extracolsep\fill}}
			\hline
			$\bm{\kappa}$ & $\bm{\theta}$  & $\bm{\sigma}$  & $\bm{\rho}$  & $\bm{r_d}$  & $\bm{r_f}$  & $\bm{T}$  & $\bm{S_0}$  & $\bm{v_0}$ \\
			\hline
			2.5 & 0.06 & 0.4 & 0.55 & 0.0198 & 0 &  1& 1&0.1683  \\
			\hline
		\end{tabular*}
	\end{table}

\begin{figure}[htb!]
\centering
\includegraphics[width=.45\columnwidth]{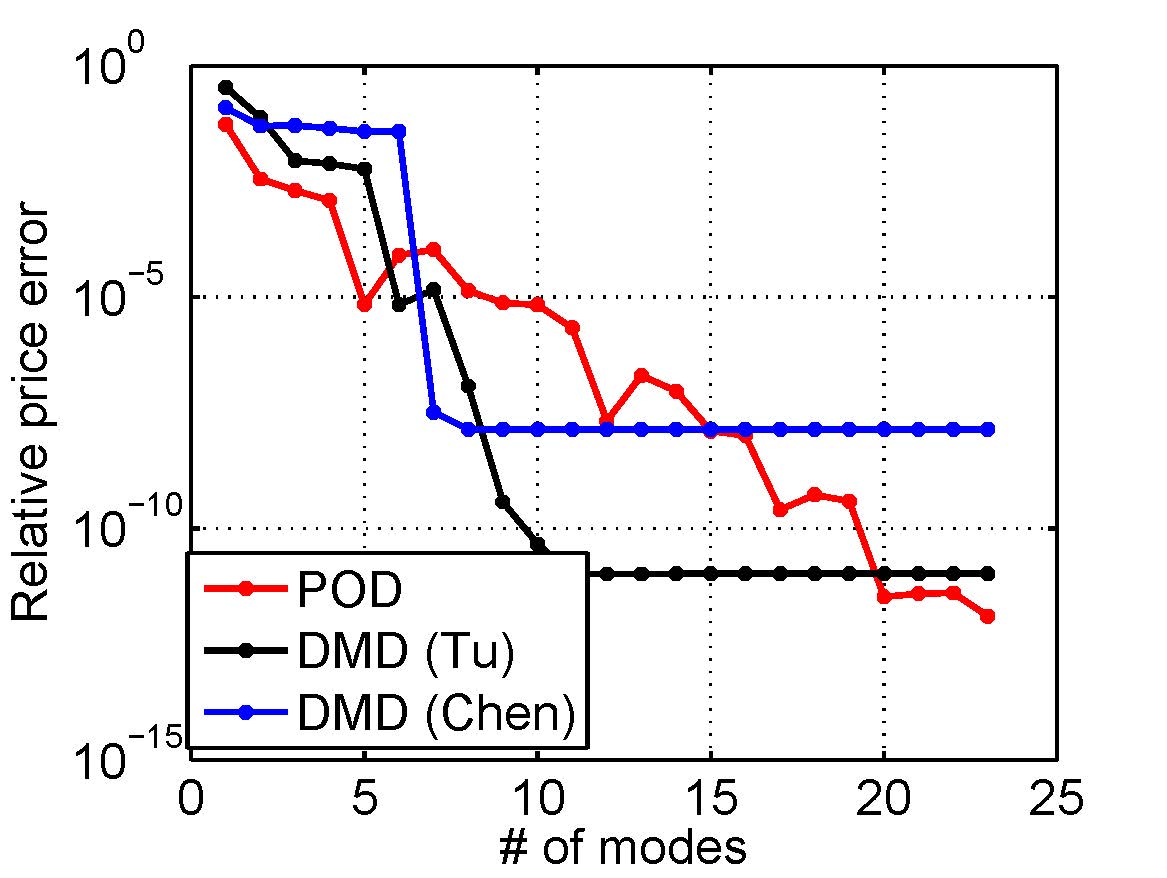}
\includegraphics[width=.45\columnwidth]{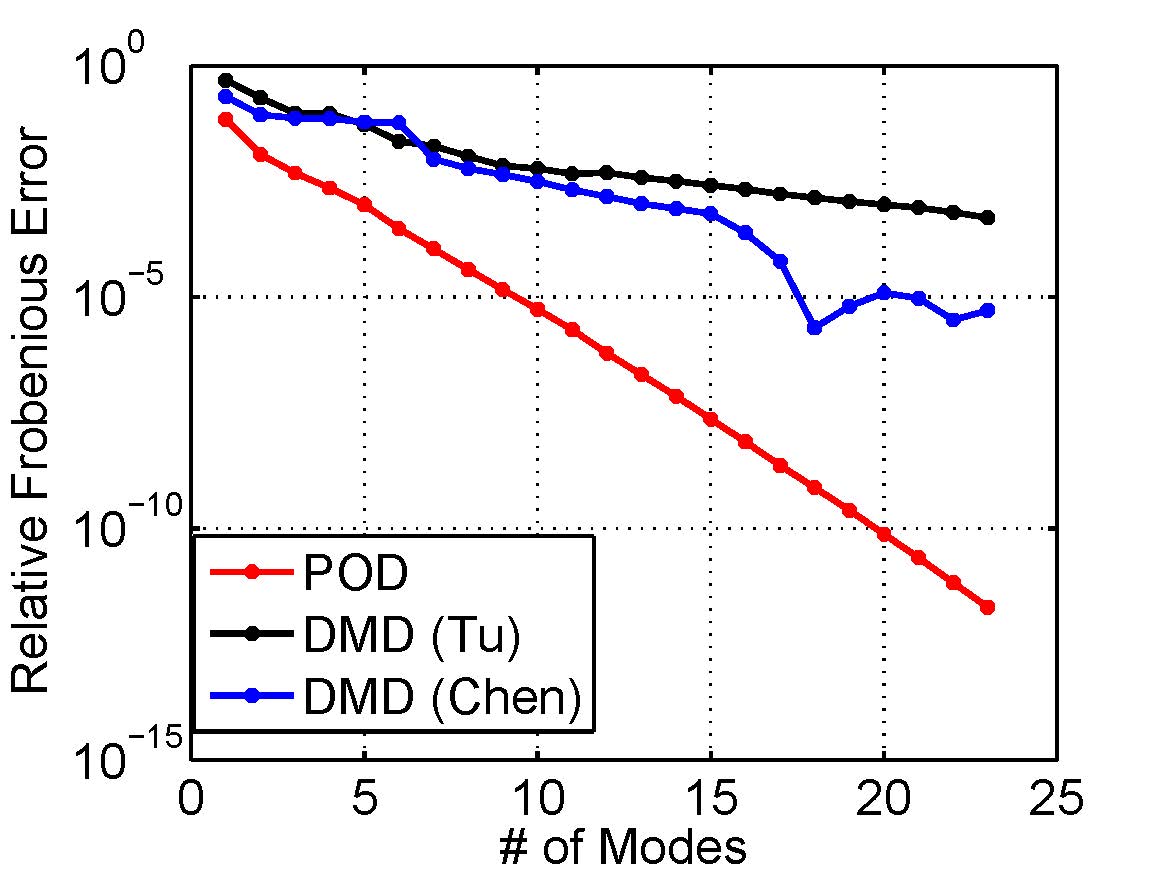}
\caption{Relative price error for $v_0=0.1683$ and $S_0=1$ (left), relative Frobenious error (right).}
\label{fig:butterfly1}
\end{figure}

\begin{figure}[htb!]
\centering
\subfloat{\includegraphics[width=.35\columnwidth]{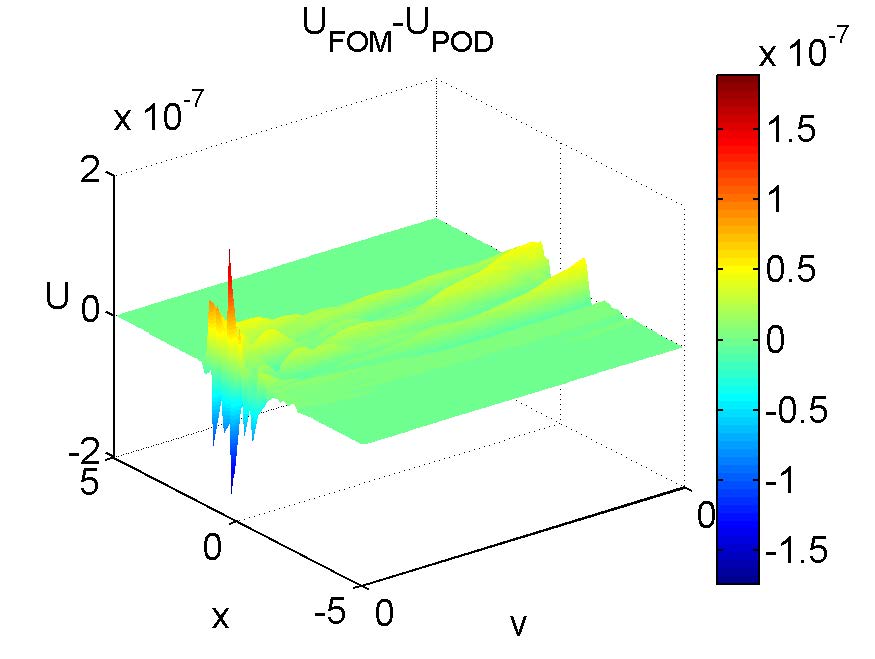}}
\subfloat{\includegraphics[width=.35\columnwidth]{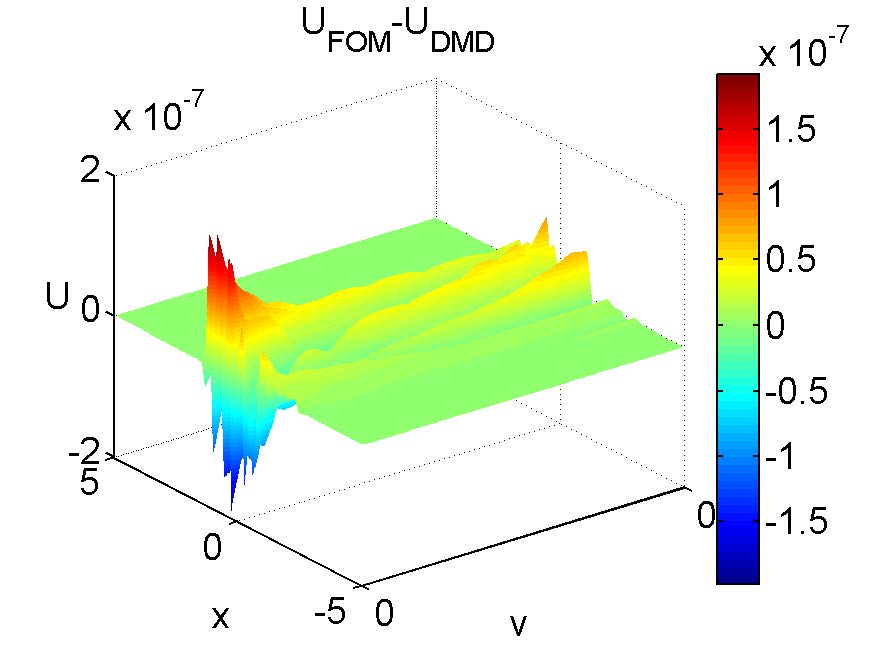}}
\subfloat{\includegraphics[width=.35\columnwidth]{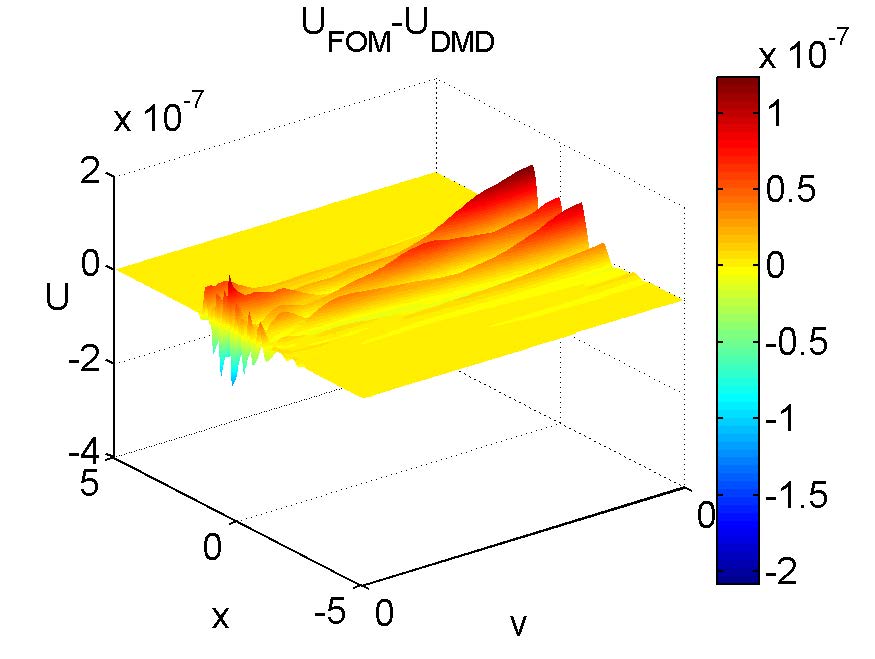}}
\caption{FOM-ROM errors: 14 POD modes (left), 23 DMD (Chen) modes (middle), 8 DMD (Tu) modes (right).}
\label{fig:butterfly2}
\end{figure}

The relative price and Frobenius errors in \figurename~\ref{fig:butterfly1} for the butterfly call option are similar as for the European call option. In this case, Tu's DMD algorithm requires less number of DMD modes than the POD at the same level of accuracy, in \figurename~\ref{fig:butterfly2}. Moreover, as described in the  European call option case, the reduced-order models could not resolve the full-order solutions in the neighborhood of $x=0$.

\subsection{Digital option}

Finally, we consider the digital call options with a discontinuous payoff \cite{england06,Lazar03}
    \begin{equation*}
    g(v,Ke^x)= \bm{1}_{\{Ke^x>K\}},
    \end{equation*}
where $K$ is the strike price of the option, which is treated as a barrier level. Precisely, if the stock price reaches the level $K$ at maturity, then the option will be worthless or it will pay 1 unit of money at time $T$.
The boundary conditions are the same as for the butterfly spread, whereas we now impose an inhomogeneous Dirichlet boundary condition at  $x=x_{\text{max}}$ ~\cite{england06}:
\begin{equation*}
u(\tau,v,x_{\text{max}})=e^{x_{\text{max}}-r_f\tau}
\end{equation*}
and an initial condition
\begin{equation*}
    u(0,v,x)= \bm{1}_{\{Ke^x>K\}}.
\end{equation*}

The discretization parameters are $\Delta t=0.01,$ $N_x=128,$ and $N_v=32$  in the domain $[0.0025,0.5]\times[-5,5].$.  The parameters of the Heston's PDE are  taken from~\cite{Winkler01}.

\begin{table}[htb!]%
		\centering
		\caption{Parameter set for the European call option.}%
		\begin{tabular*}{425pt}{@{\extracolsep\fill}cccccccccc@{\extracolsep\fill}}
			\hline
			$\bm{\kappa}$ & $\bm{\theta}$  & $\bm{\sigma}$  & $\bm{\rho}$  & $\bm{r_d}$  & $\bm{r_f}$  & $\bm{T}$  & $\bm{S_0}$  & $\bm{K}$  & $\bm{v_0}$ \\
			\hline
			2.5 & 0.06 & 0.5 & -0.1 & log(1.052) & log(1.048) & 0.25 & 1 &  1 & 0.05225  \\
			\hline
		\end{tabular*}
	\end{table}

\begin{figure}[htb!]
\centering
\includegraphics[width=.45\columnwidth]{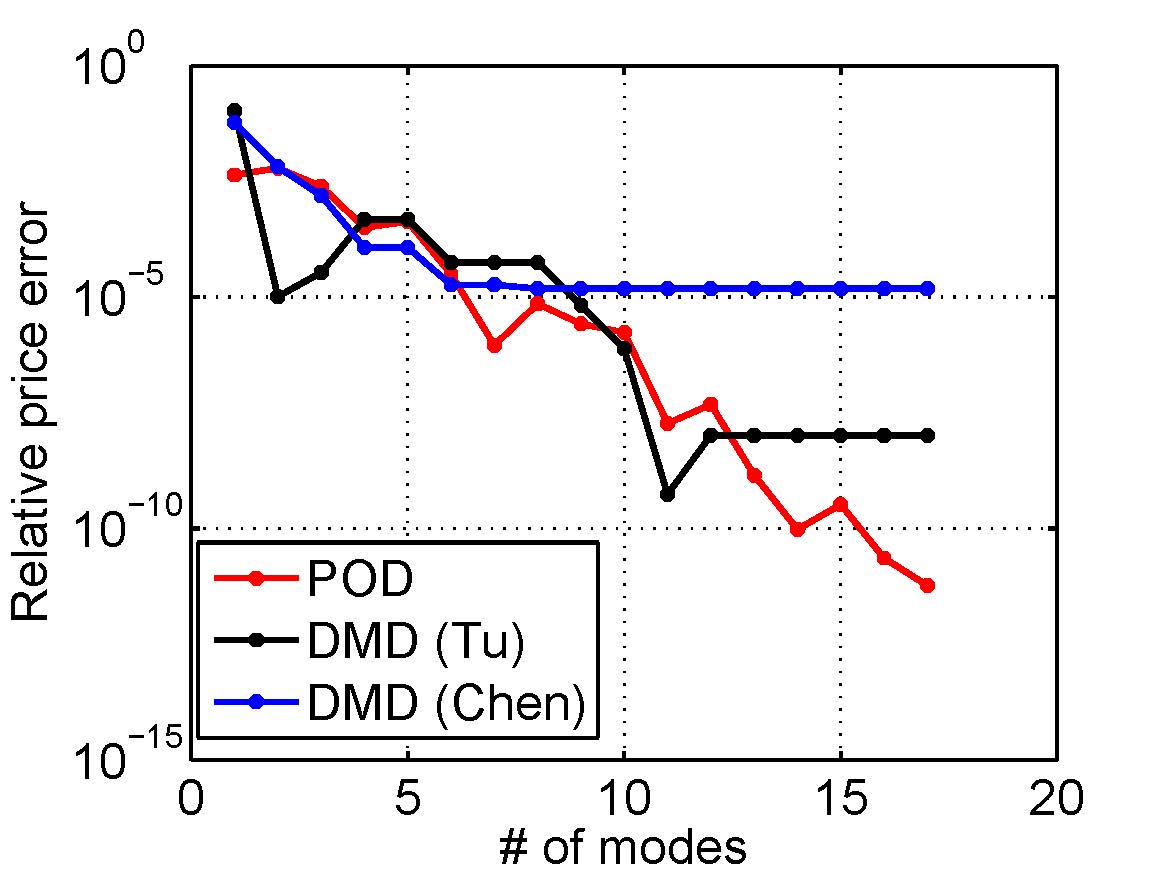}
\includegraphics[width=.45\columnwidth]{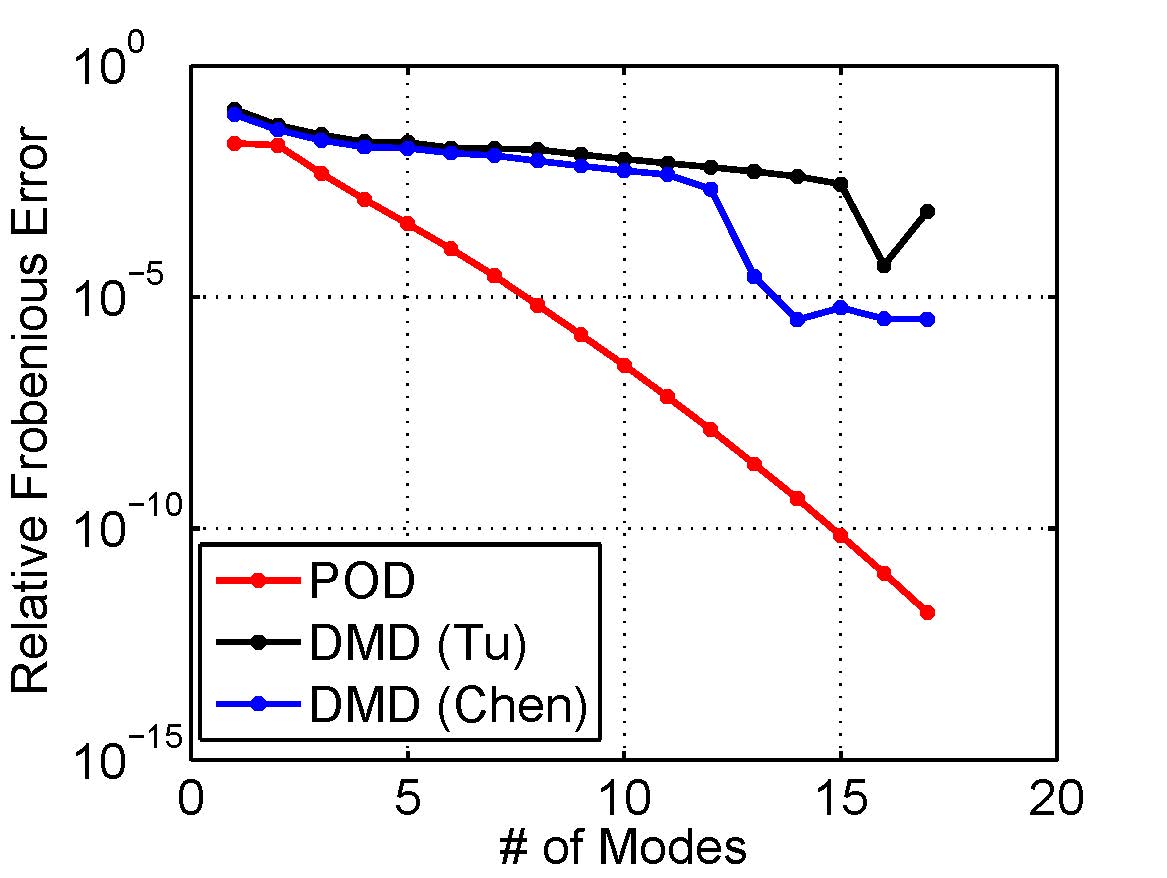}
\caption{Relative price error for $v_0=0.05225$ and $S_0=1$ (left), relative Frobenious error (right).}
\label{fig:digital1}
\end{figure}

\begin{figure}[htb!]
\centering
\subfloat{\includegraphics[width=.35\columnwidth]{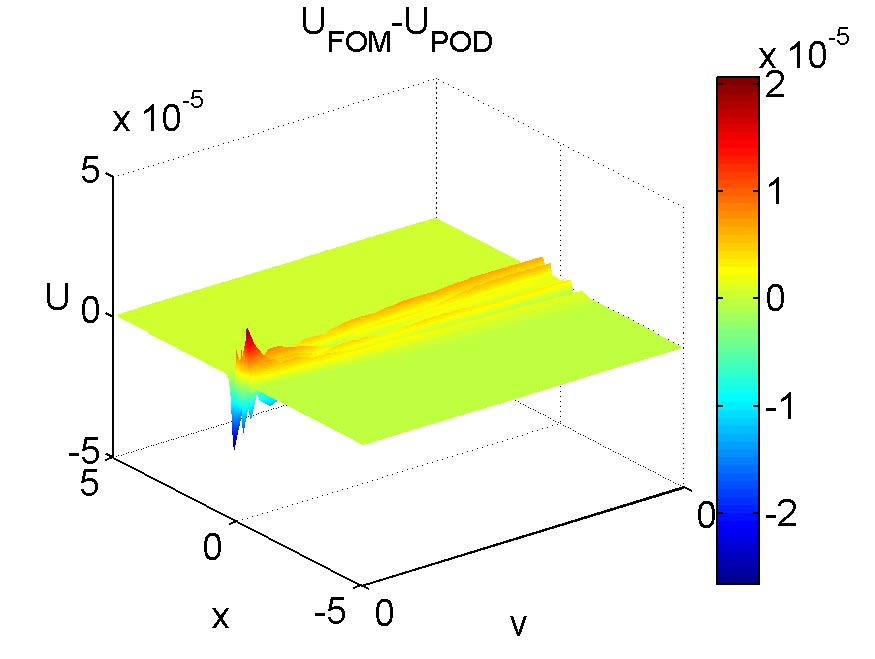}}
\subfloat{\includegraphics[width=.35\columnwidth]{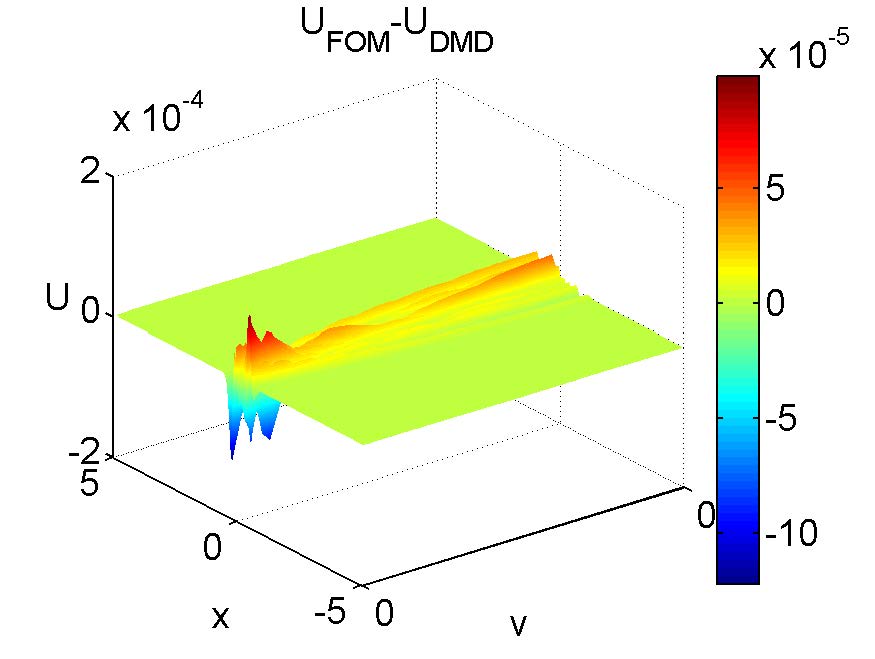}}
\subfloat{\includegraphics[width=.35\columnwidth]{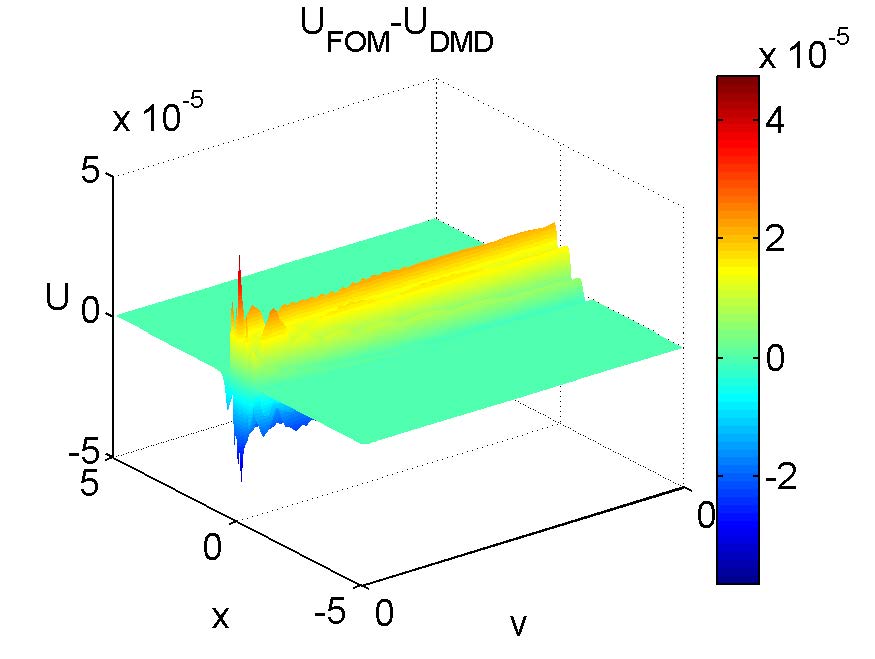}}
\caption{FOM-ROM errors: 9 POD modes (left), 17 DMD (Chen) modes (middle), 9 DMD (Tu) modes (right).}
\label{fig:digital2}
\end{figure}

The relative price and Frobenious errors and the ROM-FOM errors in \figurename~\ref{fig:digital1} and \figurename~\ref{fig:digital2} show the same behavior as for the European and butterfly spread.

\section{Conclusions}
The comparison of the POD and DMD reduced order solutions for Heston's PDE for three different options reveals that in general, the POD behaves better in terms of accuracy.  But the DMD performs better in terms of the computational cost. A selection between the two methods for a specific option should be based on balancing the accuracy of ROMs and computational cost.


\end{document}